\documentstyle[12pt]{amsart}

\author{Jie Wu}
\address{MSRI\\
1000 Centennial Drive\\
Berkeley, CA 94720}
\email{jwu@@msri.org}
\title{On Combinatorial Calculations for the James--Hopf maps}
\thanks{Research at MSRI is supported in part by NSF grant DMS-9022140.}

\newtheorem{theorem}{Theorem}[section]
\newtheorem{definition}[theorem]{Definition}

\newtheorem{notation}[theorem]{Notation}

\newtheorem{example}[theorem]{Example}

\newtheorem{lemma}[theorem]{Lemma}
\newtheorem{remark}[theorem]{Remark}

\newtheorem{proposition}[theorem]{Proposition}
\newtheorem{corollary}[theorem]{Corollary}

\begin{document}
\maketitle
\begin{abstract}
We give some formulas of the James-Hopf maps by using combinatorial
methods. An application is to give a product decomposition of the spaces
$\Omega\Sigma^2(X)$.
\end{abstract}

\section{Introduction}
In this paper, we give some formulas of the James-Hopf maps by using
combinatorial methods.\\
\par
 Let X be a
pointed space. The James-Hopf map $H_k:J(X)\rightarrow J(X^{(k)})$ is defined
by setting
$$
H_k(x_1x_2\cdots x_n)=\prod_{1\leq i_1<i_2<\cdots<i_k\leq n}(x_{i_1}x_{i_2}\cdots x_{i_k})
$$
with left lexicographical order in the product, where $X^{(k)}$ is the $k$-fold self smash product of $X$
and $J(X)$ is the James construction of $X$. The first result is
as follows.

\begin{theorem}
Let $X=\Sigma X'$ be  a suspension of a CW-complex $X'$. Then the
 composite
$$
\hspace{0.5in}
\begin{array}{cccccccccc}
&J(\vee_{l\geq k}X^{(l)}) & \stackrel{{\tilde {\vee_{l\geq k}S_l}}}{\rightarrow}&JX&\stackrel{H_k}{\rightarrow} & JX^{(k)} \\  
\end{array}
$$
is a loop map, where $S_l=[[E,\cdots,E]:X^{(l)}\rightarrow JX$ is
 the $l$-fold Samelson product and
${\tilde {\vee_{l\geq k}S_l}}:J(\vee_{l\geq k}X^{(l)})\rightarrow JX$ is the homomorphism of topological
monoids induced by $\vee_{l\geq k}S_l$.
\end{theorem}

\begin{corollary}
Let $X=\Sigma X'$ be a suspension of a CW-complex $X'$. Then 
$$
\hspace{0.5in}
\begin{array}{cccccccccccccccccc}
&JX^{(l)}&\stackrel{{\tilde S_l}}{\rightarrow}&JX&\stackrel{H_k}{\rightarrow}&JX^{(k)} \\
\end{array}
$$
is a loop map for $l\geq k$.
\end{corollary}

The second calculation is to give a decomposition of the compositions of James-Hopf
maps.

\begin{notation}
The map $L_{k.l}: X^{(kl)}\rightarrow JX^{(kl)}$ is defined by
$$
L_{k,l}(x_1\wedge x_2\wedge\cdots\wedge x_kl)=
$$
$$
\prod_{1\leq j_1<\cdots<j_l=kl,1\leq l_1^{j_s}<\cdots<l_k^{j_s}=j_s}(x_{l_1^{j_1}}\wedge\cdots\wedge x_{l_k^{j_1}}\wedge\cdots\wedge x_{l_1^{j_s}}\wedge\cdots\wedge x_{l_k^{j_l}})
$$
where $(l_1^{j_1},\cdots,l_k^{j_1},\cdots,l_1^{j_l},\cdots,l_k^{j_l})$ runs over shuffles of $(1,\cdots,kl)$ with left lexicographical oder.

Let ${\tilde L_{k,l}}:JX^{(kl)}\rightarrow J^{(kl)}$ denote the homomorphism of topological monoids induced
by $L_{k,l}$.
\end{notation}

\begin{proposition}
Let $X=\Sigma X'$ be a suspension of a CW-complex $X'$. Then
$$
H_l\circ H_k\simeq {\tilde L_{k,l}}\circ H_{kl}:JX\rightarrow JX^{(kl)}
$$
\end{proposition}

An application is to give a product decomposition of $\Omega\Sigma X$ if $X$ is a
suspension.

\begin{notation}
Let $X=\Sigma X'$. The map $\beta_n:X^{(n)}\to X^{(n)}$ is defined by induction 
$\beta_n=\beta_{n-1}\wedge1-(1,2,\dots,n)\circ(\beta_{n-1}\wedge1)$ and $\beta_2=id-(1,2)$, 
where $(1,2,\dots,n)\in\Sigma_n$ acts on $X^{(n)}$by permutating coordinates.
 Let $p$ be a prime.
 Denote $L_n(X)=hocolim_{{1\over n}\beta_n}X^{(n)}$ localized at $p$ if $n\not\equiv 0$ (mod $p$).
\end{notation}

Notice that the mod $p$ homology of the spaces $L_n(X)$ are represented by the Lie elements of
weight $n$ in the tensor algebra $T({\overline H}_*(X, {\bf Z}/p{\bf Z})$.

\begin{theorem}
Let $X=\Sigma X'$ be a suspension of a CW-complex $X'$. and let $1<k_1<k_2<\cdots$ be a sequence of
integers so that \\
(1) $k_j\not\equiv0(p)$ for each $j\geq1$\\
(2) $k_j$ is not a multiple of any $k_i$ else for each $j$.\\
Then
$$
JX\simeq \prod_jJ(L_{k_j}(X))\times?
$$
localized at the prime $p$.
\end{theorem}

The article is organized as follows. In Section 2, we introduce the groups $K_n(k)(X)$. The
combinatorial calculations and the proofs of Theorem 1.1 and Proposition 1.4 are given
in Section 3. The proof of Theorem 1.6 is given in Section 4.

\section{The groups $K_n(k)(X)$}

In this section, we consider certain subgroup in the group $[X^n,J(X^{(k)}]$, where
$X^n$ is the $n$-fold self Cartesion product of $X$ and $X^{(k)}$ is the $K$-fold self
smash product of $X$. When $k=1$, these groups have been studied by Fred Cohen [C3].

\begin{notation}
Let $X=\Sigma X'$ be the suspension of a space. Denote ${\bf K_n(k)(X)}$ the subgroup of
$[X^n,JX^{(k)}]$ generated by the homotopy classes of the maps
$$
\pi_{i_1\cdots i_k}: X^n\rightarrow JX^{(k)}
$$
with $\pi_{i_1\cdots i_k}(x_1,\cdots,x_n)=E(x_{i_1}\wedge\cdots\wedge x_{i_k})$, where $X^n$ is the $n$-fold self Cartesian
product of $X$, $X^{(k)}$ is the $k$-fold self smash product of $X$ and $E:X^{(k)}\rightarrow JX^{(k)}$
is the suspension. If there is no confusion, we simply denote $K_n(k)$ for $K_n(k)(X)$ and
denote $K_n$ or $K_n(X)$ for $K_n(1)(X)$. Denote ${\bf \{x_{i_1}|x_{i_2}|\cdots|x_{i_k}\}}$ for the homotopy class of $\pi_{i_1\cdots i_k}$
and denote ${\bf \{x^{[n_1]}_{i_1}|x^{[n_2]}_{i_2}|\cdots|x^{[n_k]}_{i_k}\}}$ for the homotopy class of the map
$f:X^n\rightarrow JX^{(k)}$ with 
$f(x_1,\cdots,x_n)=E([n_1](x_{i_1})\wedge[n_2](x_{i_2})\wedge\cdots\wedge[n_k](x_{i_k}))$, where
$[n]:X\rightarrow X$ is the composite
$$
\hspace{0.5in}
\begin{array}{cccccccccccc}
& X&\stackrel{\psi}{\rightarrow}&\bigvee^nX & \stackrel{\bigtriangledown}{\rightarrow} & X, \\  
\end{array}
$$
where $\psi$ is the comultiplication and $\bigtriangledown$ is the fold map.
\end{notation}

\begin{lemma}
Let $X=\Sigma X'$. Then, in the group $[X^n,JX^{(k)}]$, the following identities hold:\\
(1) 
$$
\{x_{i_1}|x_{i_2}|\cdots|x_{i_k}\}=1
$$
 if $i_s=i_t$ for some $1\leq s<t\leq k$\\
(2) 
$$
\{x_{i_1}^{[n_1]}|x_{i_2}^{[n_2]}|\cdots|x_{i_k}^{[n_k]}\}=\{x_{i_1}|x_{i_2}|\cdots|x_{i_k}\}^{n_1n_2\cdots n_k}
$$\\
(3) 
$$
[[\{x_{i_1}|x_{i_2}|\cdots|x_{i_k}\},\{x_{i_{k+1}}|x_{i_{k+2}}|\cdots|x_{i_{2k}}\},\cdots,\{x_{i_{(l-1)k+1}}|x_{i_{(l-1)k+2}}|\cdots|x_{i_{lk}}\}]=1
$$
if $i_s=i_t$ for some $1\leq s<t\leq kl$, where $[[a_1,a_2,\cdots,a_l]=[\cdots[a_1,a_2],\cdots,],a_l]$ 
with $[x,y]=x^{-1}y^{-1}xy$.\\
(4) 
$$
[[\{x_{i_1}|x_{i_2}|\cdots|x_{i_k}\}^{n_1},\{x_{i_{k+1}}|x_{i_{k+2}}|\cdots|x_{i_{2k}}\}^{n_2},\cdots,\{x_{i_{(l-1)k+1}}|x_{i_{(l-1)k+2}}|\cdots|x_{i_{lk}}\}^{n_l}]=
$$
$$
[[\{x_{i_1}|x_{i_2}|\cdots|x_{i_k}\},\{x_{i_{k+1}}|x_{i_{k+2}}|\cdots|x_{i_{2k}}\},\cdots,\{x_{i_{(l-1)k+1}}|x_{i_{(l-1)k+2}}|\cdots|x_{i_{lk}}\}]^{n_1n_2\cdots n_l}
$$
\end{lemma}
\noindent{\em Proof:} (1) By a shuffle map, we may assume that $i_1=i_2$. Notice that 
$\{x_{i_1}|x_{i_2}|\cdots|x_{i_k}\}$ is represented by the composite

$$
\hspace{0.5in}
\begin{array}{cccccccccccc}
& X^n&\stackrel{p_{i_1\cdots i_k}}{\rightarrow}&X^{(k)} & \stackrel{E}{\rightarrow} & JX^{(k)} \\  
\end{array}
$$
where $p_{i_1\cdots i_k}(x_1,\cdots,x_n)=x_{i_1}\wedge x_{i_2}\wedge\cdots\wedge x_{i_k}$. Since
$i_1=i_2$, there is a commutative diagram

$$
\hspace{0.5in}
\begin{array}{cccccccccc}
& X^n &\stackrel{p_{i_1i_2\cdots i_k}}{\rightarrow}&X^{(k)}                &\stackrel{E}{\rightarrow} &JX^{(k)} \\  
&\Vert&                                            & {\bar \Delta}\uparrow &                          &  \\  
& X^n&\stackrel{p_{i_1i_3\cdots i_k}}{\rightarrow} & X^{(k-1)}              &                         & \\  
\end{array}
$$
where ${\bar\Delta}:X\rightarrow X\wedge X$ is the reduced diagonal map. Since $X=\Sigma X'$ is a coH-space,
${\bar \Delta}$ is null. Thus (1) follows.\\
\\
(2) The element $\{x_{i_1}^{[n_1]}|x_{i_2}^{[n_2]}|\cdots|x_{i_k}^{[n_k]}\}$ is represented by the composite

$$
\hspace{0.5in}
\begin{array}{cccccccccc}
& X^n & \stackrel{p_{i_1\cdots i_k}}{\rightarrow}&X^{(k)}&\stackrel{[n_1]\wedge\cdots\wedge[n_k]}{\rightarrow}&X^{(k)}&\stackrel{E}{\rightarrow} & JX^{(k)} \\  
\end{array}
$$
The assertion follows from the following homotopy commutative diagram

$$
\hspace{0.5in}
\begin{array}{cccccc}
& X^{(k)} & \stackrel{E}{\rightarrow} & JX^{(k)} \\  
&\phi_1\uparrow &  & \uparrow \phi_2 \\  
&X^{(k)} &\stackrel{E}{\rightarrow} &JX^{(k)},\\  
\end{array}
$$
where $m:JX^{(k)}\rightarrow JX^{(k)}$ is the power map of degree $m$, $\phi_1=[n_1]\wedge\cdots\wedge[n_k]$
and $\phi_2=n_1\cdots n_k$.\\
\\
(3) The element $[[\{x_{i_1}|x_{i_2}|\cdots|x_{i_k}\},\cdots,\{x_{i_{(l-1)k+1}}|x_{i_{(l-1)k+2}}|\cdots|x_{i_{lk}}\}]$ is represented
by the composite

$$
\hspace{0.5in}
\begin{array}{cccccccccccc}
& X^n & \stackrel{q}{\rightarrow}(X^{(k)})^{(l)}&\stackrel{\phi}{\rightarrow} & JX^{(k)} \\  
\end{array}
$$
where $\phi=[[E, E],\cdots,E]$ and 
$$
q(x_1,x_2,\cdots,x_n)=(x_{i_1}\wedge\cdots\wedge x_{i_k})\wedge(x_{i_{k+1}}\wedge\cdots\wedge x_{i_{2k}})\wedge\cdots\wedge(x_{i_{(l-1)k+1}}\wedge\cdots\wedge x_{i_{lk}}).
$$
Since $i_s=i_t$ for some $1\leq s<t\leq kl$, there is a commutative diagram

$$
\hspace{0.5in}
\begin{array}{cccccc}
& X^n & \stackrel{q}{\rightarrow} & X^{(kl)} \\  
& \Vert &  & \uparrow {\bar\Delta} \\  
& X^n&\stackrel{q'}{\rightarrow} & X^{(kl-1)} \\  
\end{array}
$$
for a choice of reduced diagonal map ${\bar \Delta}$. Thus (3) follows.\\
\\
(4) The element  $[[\{x_{i_1}|x_{i_2}|\cdots|x_{i_k}\}^{[n_1]},\cdots,\{x_{i_{(l-1)k+1}}|x_{i_{(l-1)k+2}}|\cdots|x_{i_{lk}}\}^{[n_l]}]$ is represented
by  the composite
$$
\hspace{0.5in}
\begin{array}{cccccccccccccc}
&X^n & \stackrel{q}{\rightarrow} & (X^{(k)})^{(l)}&\stackrel{E^{(l)}}{\rightarrow}&(JX^{(k)})^{(l)}&\stackrel{n_1\wedge\cdots\wedge n_l}{\rightarrow}&(JX^{(k)})^{(l)}&\stackrel{[[id,\cdots,id]}{\rightarrow}&JX^{(k)} \\  
\end{array}
$$
where $q(x_1,\cdots,x_n)=x_{i_1}\wedge\cdots\wedge x_{i_{kl}}$. the assertion follows from the homotopy
commutative diagram

$$
\hspace{0.5in}
\begin{array}{cccccccccccc}
& (JX^{(k)})^{(l)}       & \stackrel{n_1\wedge\cdots\wedge n_l}{\rightarrow}   & (JX^{(k)})^{(l)}   &\stackrel{[[id,\cdots,id]}{\rightarrow}&JX^{(k)} \\  
&\psi \uparrow      &                                                     & \uparrow \psi &                                       &\Vert \\  
& (X^{(k)})^{(l)}        &\stackrel{[n_1]\wedge\cdots\wedge[n_l]}{\rightarrow} & (X^{(k)})^{(l)}    &\stackrel{\phi}{\rightarrow}  &JX^{(k)} \\ 
&\phi\downarrow &                                                     &                    &                                       &\Vert    \\
&JX^{(k)}                &\stackrel{=}{\rightarrow}                            &JX^{(k)}            &\stackrel{n_1n_3\cdots n_l}{\rightarrow}&JX^{(k)},\\
\end{array}
$$
where $\phi=[[E, E], \cdots, E]$ and $\psi=(E)^{(l)}$.

\begin{lemma}
Let $X=\Sigma X'$ and let $H_k:JX\rightarrow JX^{(k)}$ be the $k$-th James-Hopf map. Let
$
(H_k)_*:[X^n,JX]\rightarrow [X^n,JX^{(k)}]
$
be the function induced by $H_k$. Then, for $y=x_{i_1}^{n_1}x_{i_2}^{n_2}\cdots x_{i_l}^{n_l}\in K_n$,
\begin{eqnarray*}
(H_k)_*(y)&=&\prod_{1\leq j_1<\cdots<j_k\leq l}\{x_{i_{j_1}}|x_{i_{j_2}}|\cdots|x_{i_{j_k}}\}^{n_{j_1}n_{j_2}\cdots n_{j_k}}\\&=&\prod_{1\leq j_1<\cdots<j_k\leq l}\{x_{i_{j_1}}^{[n_{j_1}]}|x_{i_{j_2}}^{[n_{j_2}]}|\cdots|x_{i_{j_k}}^{[n_{j_k}]}\}
\end{eqnarray*}
in $K_n(k)$ with left lexicographical order.
\end{lemma}
\noindent{\em Proof:} The element $y$ is represented  by the composite

$$
\hspace{0.5in}
\begin{array}{cccccccccccc}
&X^n & \stackrel{q}{\rightarrow}&X^l&\stackrel{\prod^l_{j=1}[n_j]}{\rightarrow}&X^l&\stackrel{p}{\rightarrow}&JX,\\
  \end{array}
$$
where $p$ is the projection and $q(x_1,\cdots,x_n)=(x_{i_1},\cdots,x_{i_l})$. The assertion follows from the definition of $H_k$
and the above lemma.

\section{Combinatorial Calculations}
In this section, we give some combinatorial calculations and give proofs of Theorem 1.1
and Proposition 1.4, where Theorem 1.1 is Theorem 3.10 and Proposition 1.4 is
Proposition 3.14.

\begin{definition}
The tensor product $\otimes:[Z,JX]\times[Z,JY]\rightarrow[Z,J(X\wedge Y)]$ is defined as follows:\\
$[f]\otimes[g]$ is represented by the composite

$$
\hspace{0.5in}
\begin{array}{cccccccccccc}
& Z & \stackrel{{\bar\Delta}}{\rightarrow}&Z\wedge Z&\stackrel{f\wedge g}{\rightarrow}&JX\wedge JY&\stackrel{c}{\rightarrow} & J(X\wedge Y) \\  
\end{array}
$$
for $f:Z\rightarrow JX$ and $g:Z\rightarrow JY$, where
$
c((x_1\cdots x_n)\wedge(y_1\cdots y_m))=(x_1\wedge y_1)\cdot(x_2\wedge y_1)\cdots(x_n\wedge y_1)\cdot(x_1\wedge y_2)\cdot(x_2\wedge y_2)\cdots(x_n\wedge y_2)\cdots(x_1\wedge y_m)\cdot(x_2\wedge y_m)\cdots(x_n\wedge y_m)
$
\end{definition}
\begin{remark}
The tensor product was introduced by F.Cohen[C2]
\end{remark}
\begin{lemma}
Let $X=\Sigma X'$ and let 
$
\otimes:[X^n,JX^{(k)}]\times[X^n,JX^{(l)}]\rightarrow[X^n,JX^{(k+l)}]
$
be the tensor product. Then
$$
-\otimes\{x_{i_1}^{[n_1]}|x_{i_2}^{[n_2]}|\cdots|x_{i_l}^{[n_l]}\}:K_n(k)\rightarrow K_n(k+l)
$$
is a group homomorphism. Furthermore, 
$$
ab\otimes\{x_{i_1}^{[n_1]}|x_{i_2}^{[n_2]}|\cdots|x_{i_l}^{[n_l]}\}=ba\otimes\{x_{i_1}^{[n_1]}|x_{i_2}^{[n_2]}|\cdots|x_{i_l}^{[n_l]}\}
$$
i.e. $-\otimes\{x_{i_1}^{[n_1]}|x_{i_2}^{[n_2]}|\cdots|x_{i_l}^{n_l]}\}$ factors though the abelianization
$K_n(k)/[K_n(k),K_n(k)]$.
\end{lemma}
\noindent{\em Proof:} Let $a=\{x_{j_1}|\cdots|x_{j_k}\}$ be a generator in $K_n(k)$. Then
\begin{eqnarray*}
a^{-1}\otimes\{x_{i_1}^{[n_1]}\cdots
x_{i_l}^{[n_l]}\}&=&\{x^{[-1]}_{j_1}|x_{j_2}|\cdots|x_{j_k}\}\otimes\{x_{i_1}^{[n_1]}\cdots x_{i_l}^{[n_l]}\}\\
&=&\{x^{[-1]}_{j_1}|x_{j_2}|\cdots| x_{j_k}| x^{[n_1]}_{i_1}|\cdots|x_{i_l}^{[n_l]}\}\cr&=&(\{x_{j_1}|x_{j_2}|\cdots|x_{j_k}\}\otimes\{x_{i_1}^{[n_1]}|\cdots x_{i_l}^{[n_l]}\})^{-1}
\end{eqnarray*}
Let $a_1,a_2,\cdots,a_s$ be sequence of generators in $K_n(k)$. Then
\begin{eqnarray*}
(a_1\cdots a_s)\otimes\{x_{i_1}^{[n_1]}|x_{i_2}^{[n_2]}|\cdots|x_{i_l}^{[n_l]}\}&=&
(a_1\otimes\{x_{i_1}^{[n_1]}|x_{i_2}^{[n_2]}|\cdots|x_{i_l}^{[n_l]}\})\cdots\\
&&( a_s\otimes\{x_{i_1}^{[n_1]}|x_{i_2}^{[n_2]}|\cdots|x_{i_l}^{[n_l]}\}).
\end{eqnarray*}
Thus $-\otimes\{x_{i_1}^{[n_1]}|x_{i_2}^{[n_2]}|\cdots|x_{i_l}^{[n_l]}\}$ is a group homomorphism.
By the Lemma 2.2,
$$
ab\otimes\{x_{i_1}^{[n_1]}|x_{i_2}^{[n_2]}|\cdots|x_{i_l}^{[n_l]}\}=ba\otimes\{x_{i_1}^{[n_1]}|x_{i_2}^{[n_2]}|\cdots|x_{i_l}^{[n_l]}\}
$$
The assertion follows.

\begin{lemma}
Let $X=\Sigma X'$ and let $(H_k)_*:[X^n,JX]\rightarrow[X^n,JX^{(k)}]$ be the function induced by 
$H_k:JX\rightarrow JX^{(k)}$. Then, for $a=x_{i_1}^{n_1}\cdots x_{i_p}^{n_p}$ and 
$y=x_{i_{p+1}}^{n_{p+1}}\cdots x_{i_{p+q}}^{n_{p+q}}$ in $K_n$,
$$
(H_k)_*(a\cdot y)=(H_k)_*(a)\prod_{j=1}^q(\prod_{1\leq l_1<\cdots<l_{s-1}<l_s=j,1\leq s\leq k}(H_{k-s})_*(a)\otimes\{x_{i_{p+l_1}}^{[n_{p+l_1}]}|\cdots| x_{i_{p+l_s}}^{[n_{p+ls}]}\}),
$$
where the order of $\prod_{1\leq l_1<\cdots<l_{s-1}<l_s=j,1\leq s\leq k} $ can be chosen to any order.
\end{lemma}
\noindent{\em Proof:} By induction on $q$. If $q=0$, the assertion is trivial. Suppose that the assertion holds for $q-1$
and $y=x_{i_{p+1}}^{n_{p+1}}\cdots x_{i_{p+q}}^{n_{p+q}}$. Let
$b$ denote $a\cdot x_{i_{p+1}}^{n_{p+1}}\cdots x_{i_{p+q-1}}^{n_{p+q-1}}$. By lemma 2.3,
$$
(H_k)_*(b\cdot x_{i_{p+q}}^{n_{p+q}})=(H_k)_*(b)\cdot\prod_{1\leq l_1<\cdots<l_{k-1}\leq p+q-1}\{x_{i_{l_1}}^{[n_{l_1}]}|\cdots|x_{i_{l_{k-1}}}^{[n_{l_{k-1}}]}|x_{i_{p+q}}^{[n_{p+q}]}\}).
$$
Let 
$$
w=\prod_{1\leq l_1<\cdots<l_{k-1}\leq p+q-1}\{x_{i_{l_1}}^{[n_{l_1}]}|\cdots|x_{i_{l_{k-1}}}^{[n_{l_{k-1}}]}|x_{i_{p+q}}^{[n_{p+q}]}\}).
$$
By Lemma 2.2, the elements in the product $w$ commute each other. Thus 
$$
w=(H_{k-1})_*(b)\otimes\{x_{i_{p+q}}^{[n_{p+q}]}\}
$$
$$
=(H_{k-1})_*(a)\otimes\{x_{i_{p+q}}^{[n_{p+q}]}\}+
$$
$$
\sum^{q-1}_{j=1}\sum_{1\leq l_1<\cdots<l_{s-1}<l_s=j,1\leq s\leq k-1}(H_{k-s-1})_*(a)\otimes\{x_{i_{p+l_1}}^{[n_{p+l_1}]}|\cdots|x_{i_{p+l_s}}^{[n_{p+l_s}]}\}\otimes\{x_{i_{p+q}}^{[n_{p+q}]}\}
$$
$$
=\sum_{1\leq l_1<\cdots<l_s=q,1\leq s\leq k}(H_{k-s})_*(a)\otimes\{x_{i_{p+l_1}}^{[n_{p+l_1}]}|\cdots| x_{i_{p+ls}}^{[n_{p+l_s}]}\}.
$$
The assertion follows.

\begin{lemma}
Let $X=\Sigma X'$ and let $(H_k)_*:[X^n,JX]\rightarrow[X^n,JX^{(k)}]$ be induced by $H_k:JX\rightarrow JX^{(k)}$. 
Then
$$
(H_k)_*([[x_{i_1}^{n_1},\cdots,x_{i_m}^{n_m}])\otimes x_j=1
$$
in $[X^n,JX^{(k+1)}]$ for $m>k\geq1$.
\end{lemma}
\noindent{\em Proof:} By induction on $k$. If $k=1$, $H_1$ is the identity and the assertion follows
from Lemma 3.3. Suppose that the assertion holds for $<k$ and consider 
$(H_k)_*([[x_{i_1}^{n_1},\cdots,x_{i_m}^{n_m}])\otimes x_j$ with $m>k$. Let 
${\tilde{<x_j>}}$ denote the subgroup of $K_n(k+1)$ generated by all of the elements 
$\{x_{a_1}|\cdots|x_{a_{k+1}}\}$ so that $a_{k+1}=j$. Then, by Lemma 2.2, 
${\tilde {<x_j>}}$ is abelian.
Thus the product is the sum in ${\tilde {<x_j>}}$. Let $z_j=x_{i_j}^{n_j}$. Then
$$
(H_k)_*([[z_1,\cdots,z_m])\otimes x_j=(H_k)_*([[z_1,\cdots,z_{m-1}]^{-1}\cdot z_m^{-1}\cdot[[z_1,\cdots,z_{m-1}]\cdot z_m)\otimes x_j
$$
$$
=\sum_{s+t=k}(H_s)_*([[z_1,\cdots,z_{m-1}]^{-1})\otimes(H_t)_*([[z_1,\cdots,z_{m-1}])\otimes x_j
$$
$$
+\sum_{s+t=k-1}(H_s)_*([[z_1,\cdots,z_{m-1}]^{-1})\otimes z^{-1}_m\otimes (H_t)_*([[z_1,\cdots,z_{m-1}])\otimes x_j
$$
$$
+\sum_{s+t=k-1}(H_s)_*([[z_1,\cdots,z_{m-1}]^{-1})\otimes (H_t)_*([[z_1,\cdots,z_{m-1}])\otimes z_m\otimes x_j
$$
$$
+\sum_{s+t=k-2} (H_s)_*([[z_1,\cdots z_{m-1}]^{-1})\otimes z_m^{-1}\otimes  (H_t)_*([[z_1,\cdots,z_{m-1}])\otimes z_m \otimes x_j.
$$ 
By Lemma 2.2,
$$
(H_s)_*([[z_1,\cdots z_{m-1}]^{-1})\otimes z_m^{-1}\otimes  (H_t)_*([[z_1,\cdots,z_{m-1}])\otimes z_m \otimes x_j=1.
$$
By induction, for $1\leq t\leq k-1$,
$$
(H_s)_*([[z_1,\cdots,z_{m-1}]^{-1})\otimes(H_t)_*([[z_1,\cdots,z_{m-1}])\otimes x_j=1,
$$
$$
(H_s)_*([[z_1,\cdots,z_{m-1}]^{-1})\otimes z^{-1}_m\otimes (H_t)_*([[z_1,\cdots,z_{m-1}])\otimes x_j=1,
$$
and
$$
(H_s)_*([[z_1,\cdots,z_{m-1}]^{-1})\otimes (H_t)_*([[z_1,\cdots,z_{m-1}])\otimes z_m\otimes x_j=1
$$
Also, by induction, we have
$$
(H_{k-1})_*([[z_1,\cdots,z_{m-1}]^{-1})\otimes z^{-1}_m\otimes x_j=1
$$
and
$$
(H_{k-1})_*([[z_1,\cdots,z_{m-1}]^{-1})\otimes z_m\otimes x_j=1.
$$
Thus
$$
(H_k)_*([[z_1,\cdots,z_m])\otimes x_j=(H_k)_*([[z_1,\cdots,z_{m-1}]^{-1})\otimes x_j+(H_k)_*([[z_1,\cdots,z_{m-1}])\otimes x_j
$$
$$
=(H_k)_*([[z_1^{-1},\cdots,z_{m-1}])\otimes x_j+(H_k)_*([[z_1,\cdots,z_{m-1}])\otimes x_j
$$
 with $m-1\geq k$. Now $(H_k)_*[[z_1,\cdots,z_{m-1}]$ is represented by
$$
\hspace{0.5in}
\begin{array}{cccccccccccc}
& X^n & \stackrel{q}{\rightarrow}&X^{(m-1)}&\stackrel{[n_1]\wedge\cdots\wedge[n_{m-1}]}{\rightarrow}&X^{(m-1)}&\stackrel{[[E,\cdots,E]}{\rightarrow}& JX&\stackrel{H_k}{\rightarrow} & JX^{(k)} \\  
\end{array}
$$
and $(H_k)_*[[z_1^{-1},\cdots,z_{m-1}]$ is represented by
$$
\hspace{0.5in}
\begin{array}{cccccccccccc}
& X^n & \stackrel{q}{\rightarrow}&X^{(m-1)}&\stackrel{[-n_1]\wedge[n_2]\wedge\cdots\wedge[n_{m-1}]}{\rightarrow}&X^{(m-1)}&\stackrel{[[E,\cdots,E]}{\rightarrow}& JX&\stackrel{H_k}{\rightarrow} & JX^{(k)} \\  
\end{array}
$$
where $q(x_1,x_2,\cdots,x_n)=x_{i_1}\wedge\cdots\wedge x_{i_{m-1}}$. Thus 
$$
(H_k)_*[[z_1^{-1},\cdots,z_{m-1}]=-(H_k)_*[[z_1,\cdots,z_{m-1}]
$$
and the assertion follows.

\begin{remark}
This lemma does not hold for $m\leq k$, e.g. $(H_1)_*(x_1)\otimes x_2=\{x_1|x_2\}$ and
$(H_2)_*([x_1,x_2])\otimes x_3=\{x_1|x_2|x_3\}\cdot\{x_2|x_1|x_3\}^{-1}$.
\end{remark}

\begin{definition}
Let $G$ be a group. The {\bf lower central series} $\Gamma^kG$ is defined by induction 
$$
\Gamma^1G=G
$$
and
$$
\Gamma^kG=[\Gamma^{k-1}G,G]
$$
for $k\geq2$.
\end{definition}

\begin{theorem}
Let $X=\Sigma X'$ and let $(H_k)_*:[X^n,JX]\rightarrow [X^n,JX^{(k)}]$ be induced by $H_k:JX\rightarrow JX^{(k)}$.
Then
$$
(H_k)_*(a\cdot y)=(H_k)_*(a)\cdot(H_k)_*(y)
$$
for $a\in\Gamma^kK_n$ and $y\in K_n$, i.e. $(H_k)_*$ restricted to $K_n$ is a $\Gamma^kK_n$-map.
\end{theorem}
 
\noindent{\em Proof:} It suffices to show that
$$
(H_k)_*([[x_{i_1}^{n_1},\cdots,x_{i_m}^{n_m}]\cdot y)=(H_k)_*[[x_{i_1}^{n_1},\cdots,x_{i_m}^{n_m}]\cdot (H_k)_*(y)
$$
for $m\geq k$.
Let $y=x_{j_1}^{l_1}\cdots x_{j_t}^{n_t}$. Let $z_j=x_{i_j}^{n_j}$ and $y_s=x_{j_s}^{l_s}$.
By Lemma 3.4,
$$
(H_k)_*([[z_1,\cdots,z_m]\cdot y_1\cdots y_t)=
$$
$$
(H_k)_*([[z_1,\cdots,z_m)\prod_{j=1}^t(\prod_{1\leq i_1<\cdots<i_{s-1}<i_s=j,1\leq s\leq k}(H_{k-s})_*([[z_1,\cdots,z_m])\otimes\{y_{i_1}|\cdots|y_{i_s}\}).
$$
By Lemma 3.5,
$$
(H_{k-s}([[z_1,\cdots,z_m])\otimes\{y_{i_1}|\cdots|y_{i_s}\}=1
$$
for $1\leq s\leq k-1$. Thus
$$
(H_k)_*([[z_1,\cdots,z_m]\cdot y_1\cdots y_t)=(H_k)_*([[z_1,\cdots,z_m)\prod_{j=1}^t(\prod_{1\leq i_1<\cdots<i_{k-1}<i_k=j}\{y_{i_1}|\cdots|y_{i_k}\}).
$$
$$
=(H_k)_*([[z_1,\cdots,z_m])(H_k)_*( y_1\cdots y_t).
$$
The assertion follows.

\begin{lemma}( see also $[B1]$)
Let $X$ be a path-connected CW-complex and let $f,g:JX\rightarrow \Omega Y$ so that 
$f|_{J_nX}\simeq g|_{J_nX}$ for each $n$. Then $f\simeq g$.
\end{lemma}
\noindent{\em Proof:} There is a homotopy equivalence
$\Phi:\Sigma\vee_{n=1}^{\infty} X^{(n)}\rightarrow\Sigma JX$ so that
$\Phi_n=\Phi|_{\Sigma\vee_{j=1}^{n}X^{(j)}}:\Sigma\vee_{j=1}^nX^{(j)}\rightarrow\Sigma J_nX$ are
homotopy equivalences. Denote $f',g':\Sigma JX\rightarrow Y$ the adjoints of $f$ and $g$, respectively.
Let $i_j:\Sigma X^{(j)}\rightarrow\Sigma\vee_{n=1}^{\infty}X^{(n)}$ be the canonical inclusion.
Since $f|_{J_nX}\simeq g|_{J_nX}$ for each $n$, $f'|_{\Sigma J_n}\circ\Phi_n\simeq g'|_{\Sigma J_n}\circ\Phi_n$
for each $n$ and therefore $f'\circ\Phi\circ i_j\simeq g'\circ\Phi\circ i_j$ for each $j$.
Let
$F_j:\Sigma X^{(j)}\wedge I^+\rightarrow Y$ be a homotopy between $f'\circ\Phi\circ i_j$ and
$g'\circ\Phi\circ i_j$. Then $F=\vee_{j=1}^{\infty}F_j:(\Sigma\vee_{j=1}^{\infty}X^{(j)})\wedge I^+\rightarrow Y$ is 
a homotopy between $f'\circ \Phi$ and $g'\circ\Phi$. The assertion follows.

\begin{theorem}
Let $X=\Sigma X'$ be  a suspension of a CW-complex $X'$. Then the composite
$$
\hspace{0.5in}
\begin{array}{cccccccccc}
&J(\vee_{l\geq k}X^{(l)}) & \stackrel{{\tilde {\vee_{l\geq k}S_l}}}{\rightarrow}&JX&\stackrel{H_k}{\rightarrow} & JX^{(k)} \\  
\end{array}
$$
is a loop map, where $S_l=[[E,\cdots,E]:X^{(l)}\rightarrow JX$ is the $l$-fold Samelson product and
${\tilde {\vee_{l\geq k}S_l}}:J(\vee_{l\geq k}X^{(l)})\rightarrow JX$ is the homomorphism of topological
monoids induced by $\vee_{l\geq k}S_l$.
\end{theorem}

\noindent{\em Proof:} Let ${\tilde \beta}:J(\vee_{l\geq k}X^{(l)})\rightarrow JX^{(k)}$ be the homomorphism
of topological monoids induced by the composite
$$
\hspace{0.5in}
\begin{array}{cccccccccc}
&\vee_{l\geq k}X^{(l)} & \stackrel{ \vee_{l\geq k}S_l}{\rightarrow}&JX&\stackrel{H_k}{\rightarrow} & JX^{(k)}. \\  
\end{array}
$$
By Lemma 3.9, it suffices to show that
$$
{\tilde\beta}|_{J_n(\vee_{l\geq k}X^{(l)})}\simeq H_k\circ{\tilde{\vee_{l\geq k}S_l}}|_{J_n(\vee_{l\geq k}X^{(l)})}
$$
for each $n$. Let $j_n:J_n(\vee_{l\geq k}X^{(l)})\rightarrow J(\vee_{l\geq k}X^{(l)})$ denote the inclusion and
$p_n:(\vee_{l\geq k}X^{(l)})^n\rightarrow J_n(\vee_{l\geq k}X^{(l)})$ the projection.
Notice that 
$$
p_n^*:[J_n(\vee_{l\geq k}X^{(l)}),\Omega Y]\rightarrow [(\vee_{l\geq k}X^{(l)})^n,\Omega Y]
$$
is a monomorphism for any space $Y$. It suffice to show that
$$
{\tilde \beta}\circ j_n\circ p_n\simeq H_k\circ{\tilde{\vee_{l\geq k}S_l}}\circ j_n\circ p_n.
$$
Now, by the splitting theorem for $\Sigma(\vee_{l\geq k}X^{(l)})^n$, it suffices to show that
$$
{\tilde \beta}\circ j_n\circ p_n\circ\Phi_{l_1\cdots l_n}
\simeq H_k\circ{\tilde{\vee_{l\geq k}S_l}}\circ j_n\circ p_n\circ\Phi_{l_1\cdots l_n}
$$
for all $l_1,l_2,\cdots, l_n\geq k$, where
 $\Phi_{l_1\cdots l_n}:X^{(l_1)}\times\cdots\times X^{(l_n)}\rightarrow (\vee_{l\geq k}X^{(l)})^n$ is
 the canonical inclusion.
Let ${\tilde K_n}$ be the subgroup of $[X^{l_1+l_2+\cdots+l_n},J(\vee_{l\geq k}X^{(l)})]$ generated by
${\tilde x_1},\cdots,{\tilde x_n}$, where ${\tilde x_i}$ is represented by the composite
$$
\hspace{0.5in}
\begin{array}{cccccccccccccccccc}
&X^{l_1+\cdots+l_n}&\stackrel{q}{\rightarrow}&\prod_{j=1}^nX^{(l_j)}&\stackrel{q}{\rightarrow}&X^{(l_i)} 
&\stackrel{i}{\hookrightarrow} &(\vee_{l\geq k}X^{(l)})^n&\stackrel{j_n\circ p_n}{\rightarrow}&J(\vee_{l\geq k}X^{(l)}), \\  
\end{array}
$$
where $q$ is the projection and $i$ is the injection.
Consider the monomorphism
$$
q*:[X^{(l_1)}\times\cdots\times X^{(l_n)},J(\vee_{l\geq k}X^{(l)})]\rightarrow[X^{l_1+\cdots+l_n},J(\vee_{l\geq k}X^{(l)})].
$$
Then 
$$
q*([j_n\circ p_n\circ\Phi_{l_1\cdots l_n}])={\tilde x_1}\cdot{\tilde x_2}\cdots{\tilde x_n},
$$
where $[f]$ is the homotopy class of $f$. Now consider the homomorphism
$$
{\tilde{\vee_{l\geq k}S_l}}_*:[X^{l_1+\cdots+l_n},J(\vee_{l\geq k}X^{(l)})]\rightarrow[X^{l_1+\cdots+l_n},JX].
$$
Then 
$$
{\tilde {\vee_{l\geq k}S_l}}_*({\tilde x_i})=[[x_{l_1+\cdots+l_{i-1}+1},\cdots,x_{l_1+l_2+\cdots+l_i}]
$$
By Theorem 3.8,
$$
(H_k)_*\circ{\tilde{\vee_{l\geq k}S_l}}_*({\tilde x_1}\cdots{\tilde
x_n})=\prod_{i=1}^n(H_k)_*\circ{\tilde{\vee_{l\geq k}S_l}}_*({\tilde x_i})
={\tilde\beta}_*({\tilde x_1}{\tilde x_2}\cdots{\tilde x_n}).
$$
Thus
$$
H_k\circ{\tilde{\vee_{l\geq k}S_l}}\circ j_n\circ p_n\circ\Phi_{l_1\cdots l_n}\circ q\simeq {\tilde \beta}\circ j_n\circ p_n\circ\Phi_{l_1\cdots l_n}\circ q
$$
and therefore 
$$
H_k\circ{\tilde{\vee_{l\geq k}S_l}}\circ j_n\circ p_n\circ\Phi_{l_1\cdots l_n}\simeq {\tilde \beta}\circ j_n\circ p_n\circ\Phi_{l_1\cdots l_n}.
$$
The assertion follows.

\begin{remark} When $k=2$, this is Lemma $2.5$ in $[C2]$. The homology view of this statement was
given in $[CT]$.
\end{remark}

\begin{proposition}
Let $X=\Sigma X'$ be a suspension of a CW-complex $X'$. Then the diagram
$$
\hspace{0.5in}
\begin{array}{cccccccccc}
& JX^{(l)}\times JX & \stackrel{{\tilde S_l}\cdot 1_{JX}}{\rightarrow} & JX \\  
&\phi\downarrow &  & \downarrow H_k \\  
&JX^{(k)}\times JX^{(k)}&\stackrel{m}{\rightarrow} & JX^{(k)}\\  
\end{array}
$$
homotopy commutes for $l\geq k$, where $\phi=H_k\circ{\tilde S_l}\times H_k$ and $m$ is the multiplication.
\end{proposition}

\noindent{\em Proof:} Let ${\bar K_{n,m}}$ be the subgroup of $[X^{nl+m},JX^{(l)}\times JX]$ generated
by ${\bar x_1},\cdots,{\bar x_n}$ and ${\bar y_1},\cdots,{\bar y_m}$, where ${\bar x_i}$ and ${\bar y_j}$
are represented by $q_i,r_j:X^{nl+m}=(X^l)^n\times X^m\rightarrow JX^{(l)}\times JX$ with
$$
q_i(a_1,\cdots,a_n;b_1,\cdots,b_m)=p_l(a_i)
$$
and
$$
r_j(a_1,\cdots,a_n;b_1,\cdots,b_m)=E(b_j)
$$
for $a_s\in X^l$ and $b_t\in X$, where $p_l$ is the composite $X^l\rightarrow X^{(l)}\rightarrow JX^{(l)}$.
Let $\mu:JX\times JX\rightarrow JX$ be the mutiplication. Consider
$$
\mu\circ({\tilde S_l}\times id)_*:[X^{nl+m},JX^{(l)}\times JX]\rightarrow [X^{nl+m},JX].
$$
Then
$ \mu\circ({\tilde S_l}\times id)_*({\bar x_i})=[[x_{il+1},\cdots,x_{(i+1)l}]$ and
$\mu\circ({\tilde S_l}\times id)_*({\bar y_j})=x_{nl+j}$.
Thus
$$
\mu\circ({\tilde S_l}\times id)_*({\bar x_1}\cdots{\bar x_n}\cdot{\bar y_1}\cdots{\bar y_m})=(\prod_{i=1}^n[[x_{il+1},\cdots,x_{(i+1)l}])\cdot x_{nl+1}\cdots x_{nl+m}.
$$
By Theorem 3.8,
\begin{eqnarray*}
&&\hskip-1in H_k\circ\mu\circ({\tilde S_l}\times id)_*({\bar x_1}\cdots{\bar x_n}\cdot{\bar y_1}\cdots{\bar y_m})\\&\qquad=&(\prod_{i=1}^n(H_k)_*[[x_{il+1},\cdots,x_{(i+1)l}])\cdot(H_k)_*( x_{nl+1}\cdots x_{nl+m})\\
&\qquad=&\mu\circ((H_k\circ{\tilde S_l})\times H_k)_*({\bar x_1}\cdots{\bar x_n}\cdot{\bar y_1}\cdots{\bar y_m}).
\end{eqnarray*}
Thus
$$
H_k\circ\mu({\tilde S_l}\times id)\circ q_{n,m}\simeq\mu\circ((H_k\circ{\tilde S_l})\times H_k)\circ q_{n,m}
$$
where $q_{n,m}$ is the composite
$X^{nl+m}=(X^l)^n\times X^m\rightarrow J_n(X^{(l)})\times J_mX\hookrightarrow JX^{(l)}\times JX$. Since
$q^*_{n,m}:[J_n(X^{(l)})\times J_mX,\Omega Y]\rightarrow [(X^l)^n\times X^m,\Omega Y]$ is a monomorphism,
$$
H_k\circ\mu({\tilde S_l}\times id)|_{J_nX^{(l)}\times J_mX}\simeq\mu\circ((H_k\circ{\tilde S_l})\times H_k)|_{J_nX^{(l)}\times J_mX}.
$$
The assertion follows from Lemma 3.9.

\begin{notation}
The map $L_{k.l}: X^{(kl)}\rightarrow JX^{(kl)}$ is defined by
$$
L_{k,l}(x_1\wedge x_2\wedge\cdots\wedge x_kl)=
$$
$$
\prod_{1\leq j_1<\cdots<j_l=kl,1\leq l_1^{j_s}<\cdots<l_k^{j_s}=j_s}(x_{l_1^{j_1}}\wedge\cdots\wedge x_{l_k^{j_1}}\wedge\cdots\wedge x_{l_1^{j_s}}\wedge\cdots\wedge x_{l_k^{j_l}}
$$
where $(l_1^{j_1},\cdots,l_k^{j_1},\cdots,l_1^{j_l},\cdots,l_k^{j_l})$ runs over shuffles of $(1,\cdots,kl)$ with left lexicographical oder.

Let ${\tilde L_{k,l}}:JX^{(kl)}\rightarrow J^{(kl)}$ denote the homomorphism of topological monoids induced
by $L_{k,l}$.
\end{notation}
 
\begin{proposition}
Let $X=\Sigma X'$ be a suspension of a CW-complex $X'$. Then
$$
H_l\circ H_k\simeq {\tilde L_{k,l}}\circ H_{kl}:JX\rightarrow JX^{(kl)}
$$
\end{proposition}
\noindent{\em Proof:} Consider $K_n\subseteq [X^n,JX]$. By Lemma 3.4,
$$
(H_k)_*(x_1\cdots x_n)=\prod_{j=1}^n\sum_{1\leq l_1^j<\cdots<l_k^j=j}\{x_{l_1^j}|\cdots|x_{l_k^j}\}.
$$
By Lemma 2.2,
$$
\{x_{l_1^{j_1}}|\cdots|x_{l_k^{j_1}}|\cdots|x_{l_1^{j_l}}|\cdots |x_{l_k^{j_l}}\}=1
$$
if some $l_t^{j_s}$ repeats. Thus
$$
(H_l)_*\circ(H_k)_*(x_1\cdots x_n)=\prod_{j=1}^n\sum_{1\leq j_1<\cdots<j_l=j,1\leq l_1^{j_s}<\cdots<l_k^{j_s}=j_s}\{x_{l_1^{j_1}}|\cdots| x_{l_k^{j_1}}|\cdots |x_{l_1^{j_l}}|\cdots| x_{l_k^{j_l}}\}
$$
so that no elements in $\{l_1^{j_1},\cdots,l_k^{j_1},\cdots,l_1^{j_l},\cdots,l_k^{j_l}\}$ repeat.
Thus $$(H_l)_*\circ(H_k)_*(x_1\cdots x_n)=({\tilde L_{k,l}})_*\circ (H_{kl})_*(x_1\cdots x_n)$$ and
$$
H_l\circ H_k|_{J_nX}\circ p_n\simeq {\tilde L_{k,l}}\circ H_{kl}|_{J_nX}\circ p_n,
$$
where $p_n:X^n\rightarrow J_nX$ is the projection.
Hence
$$
H_l\circ H_k|_{J_nX}\simeq {\tilde L_{k,l}}\circ H_{kl}|_{J_nX}
$$
for each $n$. The assertion follows.

\begin{remark}
If $X=\Sigma X'$, the order of the product in the definition of $L_{k,l}$ can be chosen in  any way.
\end{remark}
\begin{example} Let $X=S^n$. Then $L_{2,2}:(S^n)^{(4)}=S^{4n}\rightarrow JS^{4n}$ is given by
$$
L_{2,2}(x_1\wedge x_2\wedge x_3\wedge x_4)=(x_2\wedge x_3\wedge x_1\wedge x_4)\cdot(x_1\wedge x_3\wedge x_2\wedge x_4)\cdot(x_1\wedge x_2\wedge x_3\wedge x_4)
$$
Thus $L_{2,2}$ is of degree $2+(-1)^n$ and the diagram
$$
\hspace{0.5in}
\begin{array}{cccccc}
& JS^n & \stackrel{H_2}{\rightarrow} & JS^{2n} \\  
& H_4\downarrow &  & \downarrow H_2 \\  
&JS^{4n} &\stackrel{J(2+(-1)^n)}{\rightarrow} & JS^{4n} \\  
\end{array}
$$
homotopy commutes. This coincides with the cohomology calculation.
\end{example}

\section{Proof of Theorem 1.6}
Let $X=\Sigma X'$ be the suspension of a space $X'$. Recall that
 The map $\beta_n:X^{(n)}\rightarrow X^{(n)}$ is defined by
induction 
$$
\beta_n=\beta_{n-1}\wedge id-(1,2,\cdots,n)\circ(\beta_{n-1}\wedge id)
$$
and
$\beta_2=id-(12)$, where $(1,2,\cdots,n) \in\Sigma_n$ acts on $X^{(n)}$ by permutation of coordinates. 
Then $\beta_n\circ\beta_n=n\beta_n$ ( see $[CW]$). If $n\not\equiv0(p)$, we denote $L_n(x)=hocolim_{n^{-1}\beta_n}X^{(n)}$ localized at $p$
for a prime $p$.
 
\noindent{\em Proof of Theorem 1.6:} Let $S_n=[[E,\cdots,E]:X^{(n)}\rightarrow JX$ be the $n$-fold Samelson product.
Then $S_n$ factors through $L_n(X)$ if $n\not\equiv0(p)$, i.e. there a homotopy commutative diagram
$$
\hspace{0.5in}
\begin{array}{cccccc}
& X^{(n)} & \stackrel{S_n}{\rightarrow} & JX \\  
&\downarrow &  & \Vert \\  
& L_n(X)&\stackrel{S_n}{\rightarrow} & JX. \\  
\end{array}
$$
Denote ${\tilde S_n}:J(L_n(X))\rightarrow JX$ the homomorphism of topological monoids induced by 
$S_n:L_n(X)\rightarrow JX$. Let $\phi$ be the composite
$$
\hspace{0.5in}
\begin{array}{cccccccccc}
& \prod^{\infty}_{j=1}J(L_{k_j}(X)) & \stackrel{\prod_j{\tilde S_{k_j}}}{\rightarrow}&\prod^{\infty}_{j=1}JX&\stackrel{multi.}{\rightarrow} & JX, \\  
\end{array}
$$
where $\prod_{j=1}^{\infty}$ is the weak infinite product. Notice that $JX^{(n)}$ is $(n-1)$-connected if 
$X$ is connected. Define $\psi:JX\rightarrow\prod^{\infty}_{j=1}J(L_{k_j}(X))$ to be the composite
$$
\hspace{0.5in}
\begin{array}{cccccc}
& JX & \stackrel{\prod_jH_{k_j}}{\rightarrow}&\prod^{\infty}_{j=1}JX^{(k_j)}&\rightarrow & \prod^{\infty}_{j=1}J(L_{k_j}(X)). \\  
\end{array}
$$
It suffices to show that 
$$\psi\circ\phi:\prod^{\infty}_{j=1}J(L_{k_j}(X))\rightarrow\prod^{\infty}_{j=1}J(L_{k_j}(X))$$
is a homotopy equivalence. Consider
$$\psi\circ\phi_*:PH_*\prod^{\infty}_{j=1}J(L_{k_j}(X))\rightarrow PH_*\prod^{\infty}_{j=1}J(L_{k_j}(X)),$$
where $H_*(-)=H_*(-;{\bf F}_p)$ the homology with coefficients in the field ${\bf F}_p$, the prime
field with $p$ elements. Notice that
$PH_*\prod^{\infty}_{j=1}J(L_{k_j}(X))\cong\bigoplus^{\infty}_{j=1}PH_*JL_{k_j}(X)$. By  Corollary 3.11,
the composite
$$
\phi_{i,j}:JX^{(k_j)}\hookrightarrow\prod^{\infty}_{j=1}JX^{(k_j)}\rightarrow JX\rightarrow\prod^{\infty}_{j=1}JX^{(k_j)}\rightarrow JX^{(k_i)}
$$
is a loop map for $k_i\leq k_j$,i.e. $i\leq j$. Since $k_j$ is not a multiple of $k_i$,
$$
\phi_{i,j}|_{H_*X^{(k_j)}}:{\bar H}_*X^{(k_j)}\rightarrow {\bar H}_*JX^{(k_i)}
$$
is zero for $k_i<k_j$ [CT, Proposition 5.3]. Now the composite
$$
\phi'_{i,i}:JL_{k_i}(X)\hookrightarrow\prod_{j=1}^{\infty}JL_{k_j}(X)\rightarrow JX\rightarrow\prod^{\infty}_{j=1}JL_{k_j}(X)\rightarrow JL_{k_i}(X)
$$
is a homotopy eqivalence [CW]. Thus the composite
$$
\hspace{-0.3in}
\begin{array}{cccccccccccccc}
&PH_*JL_{k_j}(X)& \stackrel{l_j}{\hookrightarrow}&\bigoplus^{\infty}_{j=1}PH_*JL_{k_j}(X)&\stackrel{\psi\circ\phi_*}{\rightarrow}&\bigoplus^{\infty}_{j=1}PH_*JL_{k_j}(X) & \stackrel{\pi_i}{\rightarrow}&PH_*J(L_{k_i}(X)) \\  
\end{array}
$$
is zero for $i<j$ and an isomorphism for $i=j$, where $\pi_i$ is the projection. To check that 
$$\psi\circ\phi_*:\bigoplus^{\infty}_{j=1}PH_*J(L_{k_j}(X))\rightarrow\bigoplus^{\infty}_{j=1}PH_*J(L_{k_j}(X))$$
is a monomorphism, suppose that
$$
(\psi\circ\phi)_*(a_k+a_{k+1}+\cdots)=0
$$
with $a_k\not=0$, where $a_j\in PH_*JL_{k_j}(X)$. Then
$$
\pi_k\circ(\psi\circ\phi)_*(a_k+a_{k+1}+\cdots)=\pi_k\circ(\psi\circ\phi)_*(a_k)+\sum_{j>k}\pi_k\circ(\psi\circ\phi)_*(a_j)
$$
$$
=\pi_k\circ(\psi\circ\phi)_*(a_k)=\pi_k\circ(\psi\circ\phi)_*\circ l_k(a_k)=0.
$$
Thus $a_k=0$ which is a contradiction. Thus $\psi\circ\phi_*$ is a monomorphism and therefore an isomorphism 
for any finite CW-complex. Similarly, $\psi\circ\phi$ is a rational isomorphism. Thus $\psi\circ\phi$ is
a homotopy equivalence for any finite CW-complex $X=\Sigma X'$. Notice that any CW-complex is a homotopy colimit
of finite CW-complexes. The assertion follows.

\begin{corollary}
Let $X=\Sigma X'$ be a suspension of a CW-complex $X'$ and let $q_1<q_2<\cdots $ be  all of the primes which are different
from the prime $p$. Then
$$
JX\simeq \prod^{\infty}_{j=1}JL_{q_j}(X)\times?
$$
localized at $p$.
\end{corollary}

\begin{corollary}
There exists a space $Z_{n+1}$ so that
$$
\Omega P^{n+1}(2)\simeq(\prod^{\infty}_{j=1}\Omega\Sigma L_{p_j}(P^n(2)))\times Z_{n+1}
$$
for $n\geq 3$, where $\{p_j\}$ is the set of odd primes.
\end{corollary}

\end{document}